\documentclass[12pt]{article}
\usepackage{times}
\usepackage{amsmath}
\usepackage{amssymb}
\usepackage{amsthm}
\usepackage{amsmath}
\usepackage{amsfonts}
\usepackage{amssymb}
\usepackage[all]{xy}
\usepackage{xcolor}
\usepackage[utf8]{inputenc}

\newtheorem{theorem}{Theorem}[section]
\newtheorem{lemma}[theorem]{Lemma}
\newtheorem{corollary}[theorem]{Corollary}
\newtheorem{proposition}[theorem]{Proposition}

\theoremstyle{definition}
\newtheorem{observation}[theorem]{Remark}
\newtheorem{definition}[theorem]{Definition}
\def \fim {{\hfill $\blacksquare$}}

\def\Z{\mathbb Z}

\def\0{\underline 0}

\begin{document}

\title {\bf Local topology of a deformation of a function-germ with a one-dimensional critical set   \footnote{Research partially supported by  FAPESP - Brazil, Grant 2015/25191-9 and 2017/18543-1. \newline $\quad$ {\it
Key-words: Brasselet number, Euler obstruction, Milnor fibre, Lê-Iomdin formulas}
\newline   {\it } }}

\vspace{1cm}
\author{Hellen Santana}

\date{\bf }
\maketitle

\begin{abstract}
The Brasselet number of a function $f$ with nonisolated singularities describes numerically the topological information of its generalized Milnor fibre. In this work, we consider two function-germs $f,g:(X,0)\rightarrow(\mathbb{C},0)$ such that $f$ has isolated singularity at the origin and $g$ has a stratified one-dimensional critical set. We use the Brasselet number to study the local topology a deformation $\tilde{g}$ of $g$ defined by $\tilde{g}=g+f^N,$ where $N\gg1$ and $N\in\mathbb{N}$. As an application of this study, we present a new proof of the Lê-Iomdin formula for the Brasselet number. 

\end{abstract}

\section*{Introduction}

\hspace{0,5cm} The Milnor number, defined in \cite{Milnor}, is a very useful invariant associated to a complex function $f$ with isolated singularity defined over an open neighborhood of the origin in $\mathbb{C}^N$. It gives numerical information about the local topology of the hypersurface $V(f)$  and compute the Euler characteristic of the Milnor fibre of $f$ at the origin.  

In the case where the function-germ has nonisolated singularity at the origin, the Milnor number is not well defined, but the Milnor fibre is, what led many authors (\cite{iomdin1974complex},\cite{le1980ensembles},\cite{BMPS}, \cite{DG}, \cite{massey2003numerical}) to study an extension for this number in more general settings. For example, if we consider a function with a one-dimensional critical set defined over an open subset of $\mathbb{C}^n$ and a generic linear form $l$ over $\mathbb{C}^n$, Iomdin gave an algebraic proof (Theorem 3.2), in \cite{iomdin1974complex}, of a relation between the Euler characteristic of the Milnor fibre of $f$ and the Euler characteristic of the Milnor fibre of $f+l^N, N\gg1$ and $N\in\mathbb{N},$ using properties of algebraic sets with one-dimensional critical locus. In \cite{le1980ensembles}, Lê proved (Theorem 2.2.2) this same relation in a more geometric approach and with a way to obtain the Milnor fibre of $f$ by attaching a certain number of $n$-cells to the Milnor fibre of $f|_{\{l=0\}}.$ 

In \cite{massey2003numerical}, Massey worked with a function $f$ with critical locus of higher dimension defined over a nonsingular space and defined the Lê numbers and cycles, which provides a way to numerically describe the Milnor fibre of this function with nonisolated singularity. Massey compared (Theorem II.4.5), using appropriate coordinates,  the Lê numbers of $f$ and $f+l^N,$ where $l$ is a generic linear form over $\mathbb{C}^n$  and $N\in\mathbb{N}$ is sufficiently large, obtaining a Lê-Iomdin type relation between these numbers. He also gave (Theorem II.3.3 ) a handle decomposition of the Milnor fibre of $f$, where the number of attached cells is a certain Lê number. Massey extended the concept of Lê numbers to the case of functions with nonisolated singularities defined over complex analytic spaces, introducing the Lê-Vogel cycles, and proved the Lê-Iomdin-Vogel formulas: the generalization of the Lê-Iomdin formulas in this more general sense.

The Brasselet number, defined by Dutertre and Grulha in \cite{DG}, also describes the local topological behavior of a function with nonisolated singularities defined over an arbitrarily singular analytic space: if $f:(X,0)\rightarrow(\mathbb{C},0)$ is a function-germ and $\mathcal{V}=\{\{0\},V_1,\ldots,V_q\}$ is a good stratification of $X$ relative to $f$ (see Definition \ref{good stratification}), the Brasselet number $B_{f,X}(0)$ is defined by \begin{eqnarray*}
B_{f,X}(0)=\sum^{q}_{i=1}\chi(V_i\cap f^{-1}(\delta)\cap B_{\epsilon})Eu_X(V_i).
\end{eqnarray*}

In \cite{DG}, the authors proved several formulas about the local topology of the generalized Milnor fibre of a function germ $f$ using the Brasselet number, like the Lê-Greuel type formula (Theorem 4.2 in \cite{DG}): $B_{f,X}(0)-B_{f,X^g}(0)=(-1)^{\dim_{\mathbb{C}}}n,$ where $n$ is the number of stratified Morse critical points of a Morsefication of $g|_{X\cap f^{-1}(\delta)\cap B_{\epsilon}}$ on $V_q\cap f^{-1}(\delta)\cap B_{\epsilon}.$ In \cite{DP}, Dalbelo e Pereira provided formulas to compute the Brasselet number of a function defined over a toric variety and in \cite{Daiane}, Ament, Nuño-Ballesteros, Oréfice-Okamoto and Tomazella computed the Brasselet number of a function-germ with isolated singularity at the origin and defined over an isolated determinantal variety (IDS) and the Brasselet number of finite functions defined over a reduced curve. More recently, in \cite{DH}, Dalbelo and Hartmann calculated the Brasselet number of a function-germ defined over a toric variety using combinatorical properties of the Newton polygons. In the global study of the topology of a function germ, Dutertre and Grulha defined, in \cite{DGglobal}, the global Brasselet numbers and the Brasselet numbers at infinity. In that paper, the authors compared the global Brasselet numbers of a function-germ $f$ with the global Euler obstruction of the fibres of $f,$ defined by Seade, Tib\u{a}r and Verjovsky in \cite{STVglobal}. They also related the number of critical points of a Morsefication of a polynomial function $f$ on an algebraic set $X$ to the global Brasselet numbers and the Brasselet numbers at infinity of $f.$ Therefore, the Brasselet number has been a useful tool in the study of the topology of function-germs and it will be the main object in this work.


We consider analytic function-germs $f,g:(X,0)\rightarrow(\mathbb{C},0),$ a Whitney stratification $\mathcal{W}$ of $X,$ suppose that $f$ has isolated singularity at the origin and $g$ has a one-dimensional stratified critical set. Consider the good stratification of $X$ induced by $f,$\linebreak $\mathcal{V}=\{W_{\lambda}\setminus X^f, W_{\lambda}\cap X^f\setminus\{0\},\{0\}, W_{\lambda}\in\mathcal{W}\}$ and suppose that $g$ is tractable at the origin with respect to $\mathcal{V}$ (see Definition \ref{definition tractable}). Let $\epsilon$ be sufficiently small such that the local Euler obstruction of $X^g$ is constant on $b_j\cap B_{\epsilon}$. In this case, we denote by $Eu_{X^g}(b_j)$ the local Euler obstruction of $X$ at a point of $b_j\cap B_{\epsilon}$ and by $B_{g,X\cap f^{-1}(\delta)}(b_j)$ the Brasselet number of $g|_{X\cap f^{-1}(\delta)}$ at a point of $b_j\cap B_{\epsilon}.$  For a deformation of $g$, $\tilde{g}=g+f^N, N\gg1,$ we prove (Proposition \ref{B g tilde, X^f=B f,X g tilde}) 
\begin{center}
    $B_{g,X^f}(0)=B_{\tilde{g},X^f}(0)=B_{f,X^{\tilde{g}}}(0).$
\end{center}
and for $0<|\delta|\ll\epsilon\ll1$ (Proposition \ref{Brasselet number of f over the fibre of g and g tilde}),
 \begin{center}
$B_{f,X^g}(0)-B_{f,X^{\tilde{g}}}(0)=\sum_{j=1}^{r}m_{f,b_j}(Eu_{X^g}(b_j)-B_{g,X\cap f^{-1}(\delta)}(b_j)).$
\end{center}

As an application of these results, we compare the Brasselet numbers $B_{g,X}(0)$ and $B_{\tilde{g},X}(0),$ and we obtain (Theorem \ref{Le-Iomdin formula for the Brasselet number}) a topological proof of the Lê-Iomdin formula for the Brasselet number,
\begin{eqnarray*}
    B_{\tilde{g},X}(0)=B_{g,X}(0)+ N\sum_{j=1}^{r}m_{f,b_j} Eu_{f,X\cap \tilde{g}^{-1}(\alpha^{\prime})}(b_j),
\end{eqnarray*}

\noindent where $0\ll|\alpha^{\prime}|\ll1$ is a regular value of $\tilde{g}$. This formula generalizes the Lê-Iomdin formula for the Euler characteristic of the Milnor fibre in the case of a function with isolated singularity. We note that an algebraic proof can be obtained using the description (see \cite{BMPS}) of the defect of a function-germ $f$ in therms of the Euler characteristic of vanishing cycles and the Lê-Vogel numbers associated to $f.$ 

In \cite{tibuar1998embedding}, Tib\u{a}r provided a bouquet decomposition for the Milnor fibre of $\tilde{g}$ and related it with the Milnor fibre of $g.$ As a consequence of this strong result, Tib\u{a}r gave a Lê-Iomdin formula to compare the Euler characteristics of these Milnor fibres. In the last section of this work, we apply our results to give an alternative proof for this Lê-Iomdin formula (see Proposition \ref{proposition Nicolas-Tibar}): \begin{eqnarray*}
\chi(X\cap\tilde{g}^{-1}(\alpha^{\prime})\cap 
    B_{\epsilon})=\chi(X\cap g^{-1}(\alpha)\cap B_{\epsilon})
    &+&N\sum_{j=1}^{r}m_{b_j}\big(1-\chi (F_j)\big),
\end{eqnarray*}
\noindent where $F_j=X\cap g^{-1}(\alpha)\cap H_j\cap D_{x_{t}}$ is the local Milnor fibre of $g|_{\{l=\delta\}}$ at a point of the branch $b_j$ and $H_j$ denotes the generic hyperplane $l^{-1}(\delta)$ passing through $x_t\in b_j,$ for $t\in\{i_1,\ldots,i_{k(j)}\}.$

\section{Local Euler obstruction and Euler obstruction of a function}

\hspace{0,5cm} In this section, we will see the definition of the local Euler obstruction, a singular invariant defined by MacPherson and used as one of the main tools in his proof of the Deligne-Grothendieck conjecture about the existence and uniqueness of Chern classes for singular varities. 

Let $(X,0)\subset(\mathbb{C}^n,0)$ be an equidimensional reduced complex analytic germ of dimension $d$ in a open set $U\subset\mathbb{C}^n.$ Consider a complex analytic Whitney stratification $\mathcal{V}=\{V_{\lambda}\}$ of $U$ adapted to $X$ such that $\{0\}$ is a stratum. We choose a small representative of $(X,0),$ denoted by $X,$ such that $0$ belongs to the closure of all strata. We write $X=\cup_{i=0}^{q} V_i,$ where $V_0=\{0\}$ and $V_q=X_{reg},$ where $X_{reg}$ is the regular part of $X.$ We suppose that $V_0,V_1,\ldots,V_{q-1}$ are connected and that the analytic sets $\overline{V_0},\overline{V_1},\ldots,\overline{V_q}$ are reduced. We write $d_i=dim(V_i), \ i\in\{1,\ldots,q\}.$ Note that $d_q=d.$ Let $G(d,N)$ be the Grassmannian manifold, $x\in X_{reg}$ and consider the Gauss map $\phi: X_{reg}\rightarrow U\times G(d,N)$ given by $x\mapsto(x,T_x(X_{reg})).$ 

\begin{definition}
The closure of the image of the Gauss map $\phi$ in $U\times G(d,N)$, denoted by $\tilde{X}$, is called \textbf{Nash modification} of $X$. It is a complex analytic space endowed with an analytic projection map $\nu:\tilde{X}\rightarrow X.$
\end{definition}

Consider the extension of the tautological bundle $\mathcal{T}$ over $U\times G(d,N).$ Since \linebreak$\tilde{X}\subset U\times G(d,N)$, we consider $\tilde{T}$ the restriction of $\mathcal{T}$ to $\tilde{X},$ called the \textbf{Nash bundle}, and $\pi:\tilde{T}\rightarrow\tilde{X}$ the projection of this bundle.

In this context, denoting by $\varphi$ the natural projection of $U\times G(d,N)$ at $U,$ we have the following diagram:

$$\xymatrix{
\tilde{T} \ar[d]_{\pi}\ar[r] & \mathcal{T}\ar[d] \\ 
\tilde{X}\ar[d]_{\nu}\ar[r] & U\times G(d,N)\ar[d]^{\varphi} \\ 
X\ar[r] & U\subseteq\mathbb{C}^N \\}  $$

Considering $\vert\vert z\vert\vert=\sqrt{z_1\overline{z_1}+\cdots+z_N\overline{z_N}}$, the $1$-differential form $w=d\vert\vert z\vert\vert^2$ over $\mathbb{C}^N$ defines a section in $T^{*}\mathbb{C}^N$ and its pullback $\varphi^{*}w$ is a $1$- form over $U\times G(d,N).$ Denote by $\tilde{w}$ the restriction of $\varphi^{*}w$ over $\tilde{X}$, which is a section of the dual bundle $\tilde{T}^{*}.$

Choose $\epsilon$ small enough for $\tilde{w}$ be a nonzero section over $\nu^{-1}(z), 0<\vert\vert z \vert\vert\leqslant\epsilon,$ let $B_{\epsilon}$ be the closed ball with center at the origin with radius $\epsilon$ and denote by $Obs(\tilde{T}^{*},\tilde{w})\in\mathbb{H}^{2d}(\nu^{-1}(B_{\epsilon}),\nu^{-1}(S_{\epsilon}),\Z)$ the obstruction for extending $\tilde{w}$ from $\nu^{-1}(S_{\epsilon})$ to $\nu^{-1}(B_{\epsilon})$ and $O_{\nu^{-1}(B_{\epsilon}),\nu^{-1}(S_{\epsilon})}$ the fundamental class in $\mathbb{H}_{2d}(\nu^{-1}(B_{\epsilon}),\nu^{-1}(S_{\epsilon}),\Z).$

\begin{definition}
The \textbf{local Euler obstruction} of $X$ at $0, \ Eu_X(0),$ is given by the evaluation $$Eu_X(0)=\langle Obs(\tilde{T}^{*},\tilde{w}),O_{\nu^{-1}(B_{\epsilon}),\nu^{-1}(S_{\epsilon})}\rangle.$$
\end{definition}

In \cite{BLS}, Brasselet, Lê and Seade proved a formula to make the calculation of the Euler obstruction easier.

\begin{theorem}(Theorem 3.1 of \cite{BLS})
Let $(X,0)$ and $\mathcal{V}$ be given as before, then for each generic linear form $l,$ there exists $\epsilon_0$ such that for any $\epsilon$ with $0<\epsilon<\epsilon_0$ and $\delta\neq0$ sufficiently small, the Euler obstruction of $(X,0)$ is equal to 

$$Eu_X(0)=\sum^{q}_{i=1}\chi(V_i\cap B_{\epsilon}\cap l^{-1}(\delta)).Eu_{X}(V_i),$$

\noindent where $\chi$ is the Euler characteristic, $Eu_{X}(V_i)$ is the Euler obstruction of $X$ at a point of $V_i, \ i=1,\ldots,q$ and $0<|\delta|\ll\epsilon\ll1.$
\end{theorem} 


Let us give the definition of another invariant introduced by Brasselet, Massey, Parameswaran and Seade in \cite{BMPS}. Let $f:X\rightarrow\mathbb{C}$ be a holomorphic function with isolated singularity at the origin given by the restriction of a holomorphic function $F:U\rightarrow\mathbb{C}$ and denote by $\overline{\nabla}F(x)$ the conjugate of the gradient vector field of $F$ in $x\in U,$ $$\overline{\nabla}F(x):=\left(\overline{\frac{\partial F}{\partial x_1}},\ldots, \overline{\frac{\partial F}{\partial x_n}}\right).$$

Since $f$ has an isolated singularity at the origin, for all $x\in X\setminus\{0\},$ the projection $\hat{\zeta}_i(x)$ of $\overline{\nabla}F(x)$ over $T_x(V_i(x))$ is nonzero, where $V_i(x)$ is a stratum containing $x.$ Using this projection, the authors constructed, in \cite{BMPS}, a stratified vector field over $X,$ denoted by $\overline{\nabla}f(x).$ Let $\tilde{\zeta}$ be the lifting of $\overline{\nabla}f(x)$ as a section of the Nash bundle $\tilde{T}$ over $\tilde{X}$, without singularity over $\nu^{-1}(X\cap S_{\epsilon}).$ Let $\mathcal{O}(\tilde{\zeta})\in\mathbb{H}^{2n}(\nu^{-1}(X\cap B_{\epsilon}),\nu^{-1}(X\cap S_{\epsilon}))$ be the obstruction cocycle for extending $\tilde{\zeta}$ as a nonzero section of $\tilde{T}$ inside $\nu^{-1}(X\cap B_{\epsilon}).$

\begin{definition}
The \textbf{local Euler obstruction of the function} $f, Eu_{f,X}(0)$ is the evaluation of $\mathcal{O}(\tilde{\zeta})$ on the fundamental class $[\nu^{-1}(X\cap B_{\epsilon}),\nu^{-1}(X\cap S_{\epsilon})].$
\end{definition}

The next theorem compares the Euler obstruction of a space $X$ with the Euler obstruction of function defined over $X.$

\begin{theorem}\label{Euler obstruction of a function formula}(Theorem 3.1 of \cite{BMPS})
Let $(X,0)$ and $\mathcal{V}$ be given as before and let \linebreak$f:(X,0)\rightarrow(\mathbb{C},0)$ be a function with an isolated singularity at $0.$ For $0<|\delta|\ll\epsilon\ll1,$ we have
 $$Eu_{f,X}(0)=Eu_X(0)-\sum_{i=1}^{q}\chi(V_i\cap B_{\epsilon}\cap f^{-1}(\delta)).Eu_X(V_i).$$
\end{theorem}



Let us now see a definition we will need to define a generic point of a function-germ. Let $\mathcal{V}=\{V_{\lambda}\}$ be a stratification of a reduced complex analytic space $X.$

\begin{definition}
Let $p$ be a point in a stratum $V_{\beta}$ of $\mathcal{V}.$ A \textbf{degenerate tangent plane of $\mathcal{V}$ at $p$} is an element $T$ of some Grassmanian manifold such that $T=\displaystyle\lim_{p_i\rightarrow p}T_{p_i}V_{\alpha},$ where $p_i\in V_{\alpha}$, $V_{\alpha}\neq V_{\beta}.$
\end{definition}

\begin{definition}
Let $(X,0)\subset(U,0)$ be a germ of complex analytic space in $\mathbb{C}^n$ equipped with a Whitney stratification and let $f:(X,0)\rightarrow(\mathbb{C},0)$ be an analytic function, given by the restriction of an analytic function $F:(U,0)\rightarrow(\mathbb{C},0).$ Then $0$ is said to be a \textbf{generic point}\index{holomorphic function germ!generic point of} of $f$ if the hyperplane $Ker(d_0F)$ is transverse in $\mathbb{C}^n$ to all degenerate tangent planes of the Whitney stratification at $0.$ 
\end{definition}

Now, let us see the definition of a Morsification of a function. 

\begin{definition}
Let $\mathcal{W}=\{W_0,W_1,\ldots,W_q\},$ with $0\in W_0,$ a Whitney stratification of the complex analytic space $X.$ A function $f:(X,0)\rightarrow(\mathbb{C},0)$ is said to be \textbf{Morse stratified} if $\dim W_0\geq1, f|_{W_0}: W_0\rightarrow\mathbb{C}$ has a Morse point at $0$ and $0$ is a generic point of $f$ with respect to $W_{i},$ for all $ i\neq0.$
\end{definition}

A \textbf{stratified Morsification}\index{holomorphic function germ!stratified Morsification of} of a germ of analytic function $f:(X,0)\rightarrow(\mathbb{C},0)$ is a deformation $\tilde{f}$ of $f$ such that $\tilde{f}$ is Morse stratified.

In \cite{STV}, Seade, Tib\u{a}r and Verjovsky proved that the Euler obstruction of a function $f$ is also related to the number of Morse critical points of a stratified Morsification of $f.$

\begin{proposition}(Proposition 2.3 of \cite{STV})\label{Eu_f and Morse points}
Let $f:(X,0)\rightarrow(\mathbb{C},0)$ be a germ of analytic function with isolated singularity at the origin. Then, \begin{center}
$Eu_{f,X}(0)=(-1)^dn_{reg},$
\end{center}
where $n_{reg}$ is the number of Morse points in $X_{reg}$ in a stratified Morsification of $f.$
\end{proposition}

\section{Brasselet number}

\hspace{0,5cm} In this section, we present definitions and results needed in the development of the results of this work. The main reference for this section is \cite{Ms1}.

Let $X$ be a reduced complex analytic space (not necessarily equidimensional) of dimension $d$ in an open set $U\subseteq\mathbb{C}^n$ and let $f:(X,0)\rightarrow(\mathbb{C},0)$ be an analytic map. We write $V(f)=f^{-1}(0).$ 

\begin{definition}\label{good stratification}
A \textbf{good stratification of $X$ relative to $f$} is a stratification $\mathcal{V}$ of $X$ which is adapted to $V(f)$ such that $\{V_{\lambda}\in\mathcal{V},V_{\lambda}\nsubseteq V(f)\}$ is a Whitney stratification of $X\setminus V(f)$ and such that for any pair $(V_{\lambda},V_{\gamma})$ such that $V_{\lambda}\nsubseteq V(f)$ and $V_{\gamma}\subseteq V(f),$ the $(a_f)$-Thom condition is satisfied, that is, if $p\in V_{\gamma}$ and $p_i\in V_{\lambda}$ are such that $p_i\rightarrow p$ and $T_{p_i} V(f|_{V_{\lambda}}-f|_{V_{\lambda}}(p_i))$ converges to some $\mathcal{T},$ then $T_p V_{\gamma}\subseteq\mathcal{T}.$
\end{definition}

If $f:X\rightarrow\mathbb{C}$ has a stratified isolated critical point and $\mathcal{V}$ is a Whitney stratification of $X,$ then \begin{equation}\label{induced stratification}
    \{V_{\lambda}\setminus X^f, V_{\lambda}\cap X^f\setminus\{0\},\{0\}, V_{\lambda}\in\mathcal{V}\}
\end{equation}

\noindent is a good stratification of $X$ relative to $f,$ called the good stratification induced by $f.$

\begin{definition}
The \textbf{critical locus of $f$ relative to $\mathcal{V}$}, $\Sigma_{\mathcal{V}}f,$ is given by the union \begin{center}$\Sigma_{\mathcal{V}}f=\displaystyle\bigcup_{V_{\lambda}\in\mathcal{V}}\Sigma(f|_{V_{\lambda}}).$\end{center}
\end{definition}

\begin{definition}
If $\mathcal{V}=\{V_{\lambda}\}$ is a stratification of $X,$ the \textbf{symmetric relative polar variety of $f$ and $g$ with respect to $\mathcal{V}$}, $\tilde{\Gamma}_{f,g}(\mathcal{V}),$ is the union $\cup_{\lambda}\tilde{\Gamma}_{f,g}(V_{\lambda}),$ where $\Gamma_{f,g}(V_{\lambda})$ denotes the closure in $X$ of the critical locus of $(f,g)|_{V_{\lambda}\setminus (X^f\cup X^g)},$  $X^f=X\cap \{f=0\}$ and $X^g=X\cap \{g=0\}. $ 
\end{definition}

\begin{definition}\label{definition prepolar}
Let $\mathcal{V}$ be a good stratification of $X$ relative to a function
$f:(X,0)\rightarrow(\mathbb{C},0).$ A function $g :(X, 0)\rightarrow(\mathbb{C},0)$ is \textbf{prepolar with respect to $\mathcal{V}$ at the origin} if the origin is a stratified isolated critical point, that is, $0$ is an isolated point of $\Sigma_{\mathcal{V}}g.$
\end{definition} 




\begin{definition}\label{definition tractable}
A function $g :(X, 0)\rightarrow(\mathbb{C},0)$ is \textbf{tractable at the origin with respect to a good stratification $\mathcal{V}$ of $X$ relative to $f :(X, 0)\rightarrow(\mathbb{C},0)$} if $dim_0 \ \tilde{\Gamma}^1_{f,g}(\mathcal{V})\leq1$ and, for all strata $V_{\alpha}\subseteq X^f$,
$g|_{V_{\alpha}}$ has no critical points in a neighbourhood of the origin except perhaps at the origin itself.

\end{definition}



We present now the definition of the Brasselet number. Let $f: (X,0)\rightarrow(\mathbb{C},0)$ be a complex analytic function germ and let $\mathcal{V}$ be a good stratification of $X$ relative to $f.$ We denote by $V_1,\ldots, V_q$ the strata of $\mathcal{V}$ that are not contained in $\{f=0\}$ and we assume that $V_1,\ldots, V_{q-1}$ are connected and that $V_{q}= X_{reg}\setminus \{f=0\}.$ Note that $V_q$ could be not connected.  
 
\begin{definition}
Suppose that $X$ is equidimensional. Let $\mathcal{V}$ be a good stratification of $X$ relative to $f.$ The \textbf{Brasselet number} of $f$ at the origin, $B_{f,X}(0),$ is defined by \begin{center}
$B_{f,X}(0)=\sum_{i=1}^{q}\chi(V_i\cap f^{-1}(\delta)\cap B_{\epsilon})Eu_X(V_i),$
\end{center}
where $0<|\delta|\ll\epsilon\ll1.$
\end{definition} 

\noindent\textbf{Remark:} If $V_q^i$ is a connected component of $V_{q},$ $Eu_X(V_q^i)=1.$

Notice that if $f$ has a stratified isolated singularity at the origin, then \linebreak $B_{f,X}(0)=Eu_{X}(0)-Eu_{f,X}(0)$ (see Theorem \ref{Euler obstruction of a function formula}).



\section{Local topology of a deformation of a function-germ with one-dimensional critical set}

\hspace{0,5cm} We begin this section with a discussion about the singular locus of the function $\tilde{g}=g+f^N$ and a description of the appropriate stratification with which we can compute explicitly the Brasselet numbers we will use. 

Let $f,g:(X,0)\rightarrow(\mathbb{C},0)$ be complex analytic function-germs such that $f$ has isolated singularity at the origin. Let $\mathcal{W}$ be the Whitney stratification of $X$ and $\mathcal{V}$ be the good stratification of $X$ induced by $f.$ Suppose that $\Sigma_{\mathcal{W}}g$ is one-dimensional and that $\Sigma_{\mathcal{W}}g\cap\{f=0\}=\{0\}.$

By Lemma 3.1 in \cite{Santana}, if $\mathcal{V}^f$ denote the set of strata of $\mathcal{V}$ contained in $\{f=0\}, \mathcal{V}^{\prime}=\{V_i\setminus\Sigma_{\mathcal{W}}g, V_i\cap\Sigma_{\mathcal{W}}g, V_i\in\mathcal{V}\}\cup\mathcal{V}^f$ is a good stratification of $X$ relative to $f,$ such that $ \mathcal{V}^{\prime\{g=0\}}$ is a good stratification of $X^g$ relative to $f|_{X^g},$ where \begin{center}$\mathcal{V}^{\prime\{g=0\}}=\Big\{V_{i}\cap\{g=0\} \setminus \Sigma_{\mathcal{W}} g, V_{i}\cap \Sigma_{\mathcal{W}} g, V_i\in\mathcal{V}\Big\}\cup\mathcal{V}^f\cap\{g=0\},$\end{center}
\noindent and $\mathcal{V}^f\cap\{g=0\}$ denotes the collection of strata of type $V^f\cap\{g=0\},$ with $V^f\in\mathcal{V}^f.$ In this whole section, we will use this good stratification of $X$ relative to $f.$
Suppose that $g$ is tractable at the origin with respect to $\mathcal{V}$ and let $\tilde{g}:(X,0)\rightarrow(\mathbb{C},0)$ be the function-germ given by $\tilde{g}(x)=g(x)+f^N(x), N\gg1.$

\begin{proposition}\label{g tilde has isolated singularity}
For a sufficiently large $N, \tilde{g}$ has a stratified isolated singularity at the origin with respect to the Whitney stratification $\mathcal{W}$ of $X.$ 
\end{proposition}
\noindent \textbf{Proof.}\noindent \textbf{Proof.} Let $x$ be a critical point of $\tilde{g},$ $U_x$ be a neighborhood of $x$ and $G$ and $F$ be analytic extensions of $g$ and $f$ to $U_x,$ respectively. If $V(x)$ is a stratum of $\mathcal{W}$ containing $x\neq0,$ \begin{eqnarray*}
d_x\tilde{G}|_{V(x)}=0\Leftrightarrow d_xG|_{V(x)}+N(F(x))^{N-1}d_xF|_{V(x)}=0
\end{eqnarray*}

If $d_xG|_{V(x)}=0,$ then $N(F(x))^{N-1}d_xF|_{V(x)}=0,$ hence $x\in\{F=0\}.$ \linebreak Then $x\in\Sigma_{\mathcal{W}} g\cap\{f=0\}=\{0\}.$  If $d_xG|_{V(x)}\neq0,$ we have $G\neq0.$ Since $d_x\tilde{G}|_{V(x)}=0,$ by Proposition 1.3 of \cite{Ms1}, $\tilde{G}=0,$ which implies that $F\neq0.$ On the other hand, if $d_xG|_{V(x)}\neq0,$ $d_xG|_{V(x)}=-N(F(x))^{N-1}d_xF|_{V(x)},$ and then $x\in\tilde{\Gamma}_{f,g}(V(x)).$ Suppose that $x$ is arbitrarily close to the origin. Since $f$ has isolated singularity at the origin, we can define for the stratum $V(x),$ the function $\beta:(0,\epsilon)\rightarrow\mathbb{R}, 0<\epsilon\ll1,$  \begin{eqnarray*}
\beta(u)=\inf \left\{\frac{||d_zg|_{V(x)}||}{||d_zf|_{V(x)}||}; z\in\tilde{\Gamma}_{f,g}(V(x))\cap\{|f|_{V(x)}(z)|=u, u\neq 0\}\right\},
\end{eqnarray*}
\noindent where $||.||$ denotes the operator norm, (defined, for each linear transformation $T:V\rightarrow W$ between normed vector fields, by $sup_{v\in V, ||v||=1}||T(v)||$).
Notice that, for each stratum $W_i\in\mathcal{W}, \tilde{\Gamma}_{f,g}(W_i)=\tilde{\Gamma}_{f,g}(W_i\setminus\{f=0\}).$ Since $g$ is tractable at the origin with respect to $\mathcal{V},\dim_0\tilde{\Gamma}_{f,g}(\mathcal{V})\leq1.$ Therefore,
$\dim_0\tilde{\Gamma}_{f,g}(W_i)=\dim_0\tilde{\Gamma}_{f,g}(W_i\setminus\{f=0\})\leq1.$ Hence $\tilde{\Gamma}_{f,g}(V(x))\cap\{|f|=u, u\neq 0\}$ is a finite number of points and $\beta$ is well defined.

Since the function $\beta$ is subanalytic, $\alpha(R)=\beta(1/R),$ for $R\gg1,$ is subanalytic. Then, by \cite{loi2003tame}, there exists $n_0\in\mathbb{N}$ such that $\frac{1}{\alpha(R)}< R^{n_0},$ which implies $\beta(1/R)>(1/R)^{n_0},$ that is, $\beta(u)> u^{n_0}.$ Hence, for $z\in\tilde{\Gamma}_{f,g}(V(x))\cap\{|f(z)|=u\},$ $u\ll1,$ we have
\begin{center}
$\frac{||d_zg|_{V(x)}||}{||d_zf|_{V(x)}||}\geq\beta(u)>u^{n_0},$ which implies, $||d_zg|_{V(x)}||>|f|_{V(x)}(z)|^{n_0}||d_zf|_{V(x)}||.$
\end{center}

On the other hand, since $N$ is sufficiently large, we can suppose $N>n_0.$ Since $\tilde{g}(z)=g(z)+f^N(z),$ we obtain using previous inequality that, for the critical point $x$ of $\tilde{g},$
\begin{eqnarray*}
N|f|_{V(x)}(x)|^{N-1}||d_xf|_{V(x)}||=||d_xg|_{V(x)}||>|f|_{V(x)}(x)|^{n_0}||d_xf|_{V(x)}||,
\end{eqnarray*}
\noindent which implies that $N|f|_{V(x)}(x)|^{N-1-n_0}>1.$ Since $x$ was taken sufficiently close to the origin, $f|_{V(x)}(x)$ is close to zero. Hence, $|f|_{V(x)}(x)|\ll1,$ which implies that $N-1-n_0<0.$ Therefore, $N\leq n_0,$ what is contradiction. So, there is no $x$ sufficiently close to the origin such that $d_x\tilde{g}=0.$ Therefore, $\tilde{g}$ has isolated singularity at the origin. \fim

We will now see how $\tilde{g}$ behaves with respect to the good stratification $\mathcal{V}$ of $X$ induced by $f.$

\begin{proposition}\label{g tilde is prepolar}
If $g$ is tractable at the origin with respect to the good stratification $\mathcal{V}$ of $X$ induced by $f,$ then $\tilde{g}$ is prepolar at the origin with respect to $\mathcal{V}.$
\end{proposition}

\noindent\textbf{Proof.} By Proposition \ref{g tilde has isolated singularity}, $\tilde{g}$ is prepolar at the origin with respect to $\mathcal{V}.$ So it is enough to verify that $\tilde{g}|_{V_i\cap\{f=0\}}$ is nonsingular or has isolated singularity at the origin, where $V_i$ is a stratum from the Whitney stratification $\mathcal{V}$ of $X.$ Suppose that $x\in\Sigma\tilde{g}|_{V_i\cap\{f=0\}}.$ Then $d_x\tilde{g}=d_xg+Nf(x)^{N-1}d_xf=0,$ which implies that $d_xg=0.$ But $g$ has no critical  point on $V_i\cap\{f=0\},$ since $g$ is tractable at the origin with respect to $\mathcal{V}.$ Therefore, $\tilde{g}$ is prepolar at the origin with respect to $\mathcal{V}.$\fim

\begin{corollary}\label{f is prepolar with respect to g tilde}
Let $\tilde{\mathcal{V}}$ be the good stratification of $X$ induced by $\tilde{g}.$ Then $f$ is prepolar at the origin with respect to $\tilde{\mathcal{V}}.$
\end{corollary}
\noindent\textbf{Proof.} Use Proposition \ref{g tilde is prepolar} and Lemma 6.1 of \cite{DG}.\fim

Using the previous results, we can relate the relative symmetric polar varieties $\tilde{\Gamma}_{f,\tilde{g}}(\mathcal{V})$ and $\tilde{\Gamma}_{f,g}(\mathcal{V}).$

\begin{observation}
Let us describe $\tilde{\Gamma}_{f,\tilde{g}}(\mathcal{V}).$ Let $\Sigma(\tilde{g},f)=\{x\in X; rk(d_x\tilde{g},d_x f)\leq 1\}.$ Since $f$ is prepolar at the origin with respect to the good stratification induced by $\tilde{g}, f|_{W_i\cap\{\tilde{g}=0\}}$ is nonsingular, for all $W_i\in\mathcal{W}, i\neq0.$ Also $\tilde{g}$ is prepolar at the origin with respect to the good stratification induced by $f,$ which implies that $\tilde{g}|_{W_i\cap\{f=0\}}$ is nonsingular, for all $W_i\in\mathcal{W}, i\neq0.$ Nevertheless, since $f$ and $\tilde{g}$ have stratified isolated singularity at the origin, $\Sigma_{\mathcal{W}}\tilde{g}\cup\Sigma_{\mathcal{W}} f=\{0\}.$ Therefore, the map $(f,g)$ has no singularities in $\{g=0\}$ or in $\{f=0\}.$ Hence, $\Sigma(\tilde{g},f)=\tilde{\Gamma}_{f,\tilde{g}}(\mathcal{V}).$ So, it is sufficient to describe $\Sigma(\tilde{g},f).$ Let $x\in\Sigma(\tilde{g},f),$ then
\vspace{-0,1cm}
\begin{eqnarray}
rk(d_x\tilde{g},d_xf)\leq1 &\Leftrightarrow & \big(d_x\tilde{g}=0\big) \ or \ \big(d_xf=0\big) \ or \ \big(d_x\tilde{g}=\lambda d_xf\big)\nonumber\\
&\Leftrightarrow & \big(d_x\tilde{g}=0\big) \ or \ \big(d_xf=0\big) \ or \  \big(d_xg=(-Nf(x)^{N-1}+\lambda)d_xf\big)\nonumber
\end{eqnarray}

Since $x\not\in\{f=0\}, d_xf\neq0.$ And since $\tilde{g}$ has isolated singularity at the origin, $d_x\tilde{g}\neq0.$ If $-Nf(x)^{N-1}+\lambda=0,$ then $d_xg=0,$ that is, $x\in\Sigma_{\mathcal{W}} g.$ If $-Nf(x)^{N-1}+\lambda\neq0,$ then $d_xg$ is a nonzero multiple of $d_xf,$ that is, $x\in\tilde{\Gamma}_{f,g}(\mathcal{V}).$ Therefore, \begin{eqnarray*}
\Sigma(\tilde{g},f)\subseteq \Sigma_{\mathcal{W}} g\cup\tilde{\Gamma}_{f,g}(\mathcal{V}).
\end{eqnarray*}

On the other hand, if $x\in\Sigma_{\mathcal{W}} g,$ then $d_xg=0,$ and \begin{eqnarray*}
d_x\tilde{g}=d_xg+Nf(x)^{N-1}d_xf=Nf(x)^{N-1}d_xf.
\end{eqnarray*}
So, $x\in\Sigma(\tilde{g},f).$ If $x\in\tilde{\Gamma}_{f,g}(\mathcal{V}), d_xg=\lambda d_xf,$ and \begin{eqnarray*}
d_x\tilde{g}=d_xg+Nf(x)^{N-1}d_xf=(\lambda+N)f(x)^{N-1}d_xf,
\end{eqnarray*}
\noindent which implies $x\in\Sigma(\tilde{g},f).$ Therefore, $\tilde{\Gamma}_{f,\tilde{g}}(\mathcal{V})=\Sigma(\tilde{g},f)=\Sigma_{\mathcal{W}} g\cup\tilde{\Gamma}_{f,g}(\mathcal{V}).$
\end{observation}

\begin{proposition}\label{B g tilde, X^f=B f,X g tilde}
Suppose that $g$ is tractable at the origin with respect to the good stratification $\mathcal{V}$ of $X$ induced by $f.$ Then, for $N\gg1,$ \begin{center}
    $B_{g,X^f}(0)=B_{\tilde{g},X^f}(0)=B_{f,X^{\tilde{g}}}(0).$
\end{center}
\end{proposition}
\noindent\textbf{Proof.} Since $\tilde{g}=g+f^N,$ over $\{f=0\}, \tilde{g}=g.$ Therefore, $B_{g,X^f}(0)=B_{\tilde{g},X^f}(0).$ On the other hand, by Corollary \ref{f is prepolar with respect to g tilde}, $f$ is prepolar at the origin with respect to the good stratification $\tilde{\mathcal{V}}$ of $X$ induced by $\tilde{g}$ and so is $\tilde{g}$ with respect to $\mathcal{V},$ by Proposition \ref{g tilde is prepolar}. Hence, by Corollary 6.3 of \cite{DG}, $B_{f,X^{\tilde{g}}}(0)=B_{\tilde{g},X^f}(0).$ \fim

\begin{corollary}\label{Bg tilda,X H=Bg,X H}
Let $l$ be a generic linear form in $\mathbb{C}^n$ and denote $l^{-1}(0)$ by $H.$ Then \begin{equation*}
    B_{g,X\cap H}(0)=B_{\tilde{g},X\cap H}(0)=Eu_{X^{\tilde{g}}}(0).
\end{equation*}
\end{corollary}
\noindent\textbf{Proof.} By \cite{Santana}, $g$ is tractable at the origin with respect to the good stratification $\mathcal{V}$ of $X$ induced by $l.$ Hence, the formula follows directly by Proposition \ref{B g tilde, X^f=B f,X g tilde}, using the equality $ B_{g,X^l}(0)=B_{\tilde{g},X^l}(0),$ and Corollary 6.6 of \cite{DG}.\fim

\begin{corollary}\label{inequalities for Euler obstruction g tilde} Let $N\in\mathbb{N}$ be a sufficiently large number.
\begin{enumerate}
\item {If $d$ is even, $Eu_{X^{\tilde{g}}}(0)\geq Eu_{X^g}(0);$}

\item {If $d$ is odd, $Eu_{X^{\tilde{g}}}(0)\leq Eu_{X^g}(0).$}
\end{enumerate}
\end{corollary}
\noindent\textbf{Proof.} Use Corollary 4.11 of \cite{Santana} and Corollary \ref{Bg tilda,X H=Bg,X H}.\fim

In order to compare the Brasselet numbers $B_{f,X^g}(0)$ and $B_{f,X^{\tilde{g}}}(0)$ we need to understand the stratified critical set of $g.$ We use the description presented in \cite{Santana}. Consider a decomposition of $\Sigma_{\mathcal{W}}g$ into branches $b_j,$
\begin{eqnarray*}
\Sigma_{\mathcal{W}} g=\bigcup_{\alpha=1}^q\Sigma g|_{W_{\alpha}}\cup\{0\}=b_1\cup\ldots\cup b_r,
\end{eqnarray*}

\noindent where $b_j\subseteq W_{\alpha},$ for some $\alpha\in\{1,\ldots,q\}.$ Notice that a stratum $W_{\alpha}$ can contain no branch and that a stratum  $V_j$ can contain more than one branch, but a branch can not be contained in two different strata. Let $\delta$ be a regular value of $f, 0<|\delta|\ll1,$ and let us write, for each $j\in\{1,\ldots,r\}, f^{-1}(\delta)\cap b_j=\{x_{i_1},\ldots,x_{i_{k(j)}}\}.$ So, in this case, the local degree $m_{f,b_j}$ of $f|_{b_j}$ is $k(j).$ Let $\epsilon$ be sufficiently small such that the local Euler obstruction of $X$ and of $X^g$ are constant on $b_j\cap B_{\epsilon}$. Denote by $Eu_{X}(b_j)$ (respectively, $Eu_{X^g}(b_j)$) the local Euler obstruction of $X$ (respectively, $X^g$) at a point of $b_j\cap B_{\epsilon}.$

\begin{observation}\label{constant over the branch}
If $\epsilon$ is sufficiently small and $x_l\in b_j, l\in\{i_1,\ldots,i_{k(j)}\},$ $B_{g,X\cap f^{-1}(\delta)}(x_l)$ is constant on $b_j\cap B_{\epsilon}$ (see Remark 4.5 of \cite{Santana}). Then we denote $B_{g,X\cap f^{-1}(\delta)}(x_l)$ by $B_{g,X\cap f^{-1}(\delta)}(b_j).$ Since $B_{g,X\cap f^{-1}(\delta)}(x_l)=Eu_{X\cap f^{-1}(\delta)}(x_l)-Eu_{g,X\cap f^{-1}(\delta)}(x_l),$ we also denote $Eu_{g,X\cap f^{-1}(\delta)}(x_l)$ by $Eu_{g,X\cap f^{-1}(\delta)}(b_j).$
\end{observation}

\begin{proposition}\label{Brasselet number of f over the fibre of g and g tilde}
Suppose that $g$ is tractable at the origin with respect to the good stratification $\mathcal{V}$ of $X$ relative to $f.$ Then, for $0<|\delta|\ll\epsilon\ll1,$ \begin{center}
$B_{f,X^g}(0)-B_{f,X^{\tilde{g}}}(0)=\sum_{j=1}^{r}m_{f,b_j}(Eu_{X^g}(b_j)-B_{g,X\cap f^{-1}(\delta)}(b_j)).$
\end{center}
\end{proposition}

\noindent\textbf{Proof.} Use Corollary 4.8 of \cite{Santana} and Proposition \ref{B g tilde, X^f=B f,X g tilde}.\fim

\begin{corollary}\label{Euler obstruction of g and g tilde}
For $0<|\delta|\ll\epsilon\ll1,$ \begin{eqnarray}
Eu_{X^g}(0)-Eu_{X^{\tilde{g}}}(0)=\sum_{j=1}^{r}m_{b_j}(Eu_{X^g}(b_j)-B_{g,X\cap l^{-1}(\delta)}(b_j)).\label{Eu g and Eu g tilde}
\end{eqnarray}
\end{corollary}
\noindent\textbf{Proof.} By \cite{Santana}, $g$ is tractable at the origin with respect to the good stratification $\mathcal{V}$ of $X$ induced by a generic linear form $l.$ Hence, the formula follows directly from Proposition \ref{Brasselet number of f over the fibre of g and g tilde}, using that $B_{l,X^g}(0)=Eu_{X^g}(0)$ and that $B_{l,X^{\tilde{g}}}(0)=Eu_{X^{\tilde{g}}}(0).$ \fim

\begin{observation}
Since $l$ is a generic linear form over $\mathbb{C}^n, l^{-1}(\delta)$ intersects $X\cap\{g=0\}$ tranversely and using Corollary 6.6 of \cite{DG}, we have $Eu_{X^g}(b_j)=Eu_{X^g\cap l^{-1}(\delta)}(b_j\cap l^{-1}(\delta))=B_{g,X\cap l^{-1}(\delta)\cap L}(b_j\cap l^{-1}(\delta)),$ where $L$ is a generic hyperplane in $\mathbb{C}^n$ passing through $x_l\in b_j\cap l^{-1}(\delta), j\in\{1,\ldots,r\}$ and $ l\in\{i_1,\ldots, i_{k(j)}\}.$ Denoting $B_{g,X\cap l^{-1}(\delta)\cap L}(b_j\cap l^{-1}(\delta))$ by $B^{\prime}_{g,X\cap l^{-1}(\delta)}(b_j),$ the formula obtained in Corollary \ref{Euler obstruction of g and g tilde} can be written as \begin{center}
$Eu_{X^g}(0)-Eu_{X^{\tilde{g}}}(0)=\sum_{j=1}^{r}m_{b_j}(B^{\prime}_{g,X\cap l^{-1}(\delta)}(b_j)-B_{g,X\cap l^{-1}(\delta)}(b_j)).$
\end{center}
\end{observation}

Let $m$ be the number of stratified Morse points of a partial Morsefication of $g|_{X\cap f^{-1}(\delta)\cap B_{\epsilon}}$ appearing on $X_{reg}\cap f^{-1}(\delta)\cap\{g\neq0\}\cap B_{\epsilon}$ and $\tilde{m}$ be the number of stratified Morse points of a Morsefication of $\tilde{g}|_{X\cap f^{-1}(\delta)\cap B_{\epsilon}}$ appearing on $X_{reg}\cap f^{-1}(\delta)\cap\{\tilde{g}\neq0\}\cap B_{\epsilon}.$ The next lemma shows how to compare $m$ and $\tilde{m}.$ In the following we keep the same description of $\Sigma_{\mathcal{W}}g$. 

\begin{corollary}\label{Comparing m and m tilde}
Suppose that $g$ is tractable at the origin with respect to the good stratification $\mathcal{V}$ of $X$ relative to $f.$ Then
\begin{eqnarray*}
    \tilde{m}=(-1)^{d-1}\sum_{j=1}^{r}m_{f,b_j} Eu_{g,X\cap f^{-1}(\delta)}(b_j)+m.
\end{eqnarray*}
\end{corollary}

\noindent\textbf{Proof.} By Theorem 3.2 of \cite{Santana}, \begin{center}
$B_{f,X}(0)-B_{f,X^g}(0)-\sum_{j=1}^{r}m_{f,b_j}(Eu_{X}(b_j)-Eu_{X^g}(b_j))=(-1)^{d-1}m,$ 
\end{center}
\noindent and by Proposition \ref{g tilde is prepolar}, $\tilde{g}$ is prepolar at the origin with respect to $\mathcal{V},$ by Theorem 4.4 of \cite{DG}, \begin{center}
$B_{f,X}(0)-B_{f,X^{\tilde{g}}}(0)=(-1)^{d-1}\tilde{m}.$
\end{center}

Using Proposition \ref{Brasselet number of f over the fibre of g and g tilde}, we obtain that  
\begin{eqnarray*}
\tilde{m}= m+(-1)^{d-1}\sum_{j=1}^{r}m_{f,b_j}(Eu_{X}(b_j)-B_{g,X\cap f^{-1}(\delta)}(b_j)).
\end{eqnarray*}

Since $f$ has isolated singularity at the origin, $f^{-1}(\delta)$ intersects each stratum out of \linebreak$\{f=0\}$ transversely. So, $Eu_X(V_i)=Eu_{X\cap f^{-1}(\delta)}(S),$ for each connected component of $V_i\cap f^{-1}(\delta).$ In particular, $Eu_X(b_j)=Eu_{X\cap f^{-1}(\delta)}(b_j).$ The formula holds by Theorem \ref{Euler obstruction of a function formula},
$Eu_X(b_j)-B_{g,X\cap f^{-1}(\delta)}(b_j)=Eu_{g,X\cap f^{-1}(\delta)}(b_j).$\fim

\begin{proposition}\label{comparing brasselet numbers}
Let $\tilde{\alpha}$ be a regular value of $\tilde{g}$ and $\alpha_t$ a regular value of $f,$ $0\ll|\tilde{\alpha}|\ll |\alpha_t|\ll1.$  If $g$ is tractable at the origin with respect to $\mathcal{V}$ relative to $f,$ then    $B_{g,X\cap f^{-1}(\alpha_t)}(b_j)=B_{f,X\cap \tilde{g}^{-1}(\tilde{\alpha})}(b_j).$
\end{proposition}
\noindent\textbf{Proof.} Let $x_t\in\{f=\alpha_t\}\cap b_j,$ $D_{x_t}$ be the closed ball with center at $x_{t}$ and radius $r_{t}, \linebreak0<|\alpha-\delta|\ll|\alpha_t|\ll r_{t}\ll1.$ We have
\begin{eqnarray*}
B_{g,X\cap f^{-1}(\alpha_t)}(x_t)&=&\sum\chi(W_i\cap f^{-1}(\alpha_t)\cap g^{-1}(\alpha-\delta)\cap D_{x_t})Eu_{X\cap f^{-1}(\alpha_t)}(W_i\cap f^{-1}(\alpha_t)) \\
&=&\sum\chi(W_i\cap f^{-1}(\alpha_t)\cap g^{-1}(\alpha-\delta)\cap D_{x_t})Eu_{X}(W_i).
\end{eqnarray*}

Let $g(x_t)=\alpha,\tilde{g}(x_t)=\alpha^{\prime}$ and $f(x_t)=\alpha_t.$ Then \begin{eqnarray*}
p\in f^{-1}(\alpha_t)\cap g^{-1}(\alpha-\delta)&\Leftrightarrow& g(p)=\alpha -\delta\ and \ f(p)=\alpha_t\\
&\Leftrightarrow& g(p)=g(x_t)-\delta \ and \ f(p)=\alpha_t\\
&\Leftrightarrow& g(p)+\alpha_t^N=\alpha+\alpha_t^N-\delta \ and \ f(p)=\alpha_t\\
&\Leftrightarrow& g(p)+f^N(p)=g(x_t)+f^N(x_t)-\delta \ and \ f(p)=\alpha_t\\
&\Leftrightarrow& \tilde{g}(p)=\tilde{g}(x_t)-\delta \ and \ f(p)=\alpha_t\\
&\Leftrightarrow& \tilde{g}(p)=\alpha^{\prime}-\delta \ and \ f(p)=\alpha_t.
\end{eqnarray*}

Therefore, denoting $\tilde{\alpha}=\alpha^{\prime}-\delta,$
\begin{eqnarray*}
B_{g,X\cap f^{-1}(\alpha_t)}(x_t)&=&\sum\chi(W_i\cap f^{-1}(\alpha_t)\cap g^{-1}(\alpha-\delta)\cap D_{x_t})Eu_{X}(W_i)\\
&=&\sum\chi(W_i\cap f^{-1}(\alpha_t)\cap \tilde{g}^{-1}(\tilde{\alpha})\cap D_{x_t})Eu_{X\cap  \tilde{g}^{-1}(\tilde{\alpha})}(W_i\cap \tilde{g}^{-1}(\tilde{\alpha})) \\
&=&B_{f,X\cap \tilde{g}^{-1}(\tilde{\alpha})}(x_t).
\end{eqnarray*}\fim

An immediate consequence of the last proposition is the following.

\begin{corollary}\label{comparing euler obstruction of funcions}
Let $\tilde{\alpha}$ be a regular value of $\tilde{g}$ and $\alpha_t$ a regular value of $f,$ $0\ll|\tilde{\alpha}|\ll |\alpha_t|\ll1.$  If $g$ is tractable at the origin with respect to $\mathcal{V},$ then    $Eu_{g,X\cap f^{-1}(\alpha_t)}(b_j)=Eu_{f,X\cap \tilde{g}^{-1}(\tilde{\alpha})}(b_j).$
\end{corollary}
\noindent\textbf{Proof.} Let $x_t\in\{f=\alpha_t\}\cap b_j,$ $D_{x_t}$ be the closed ball with center at $x_{t}$ and radius $r_{t}, \linebreak0<|\alpha^{\prime}|\ll|\alpha_t|\ll r_{t}\ll1.$ The equality holds by Proposition \ref{comparing brasselet numbers}.\fim

\section{Lê-Iomdin formula for the Brasselet number}

\hspace{0,5cm} Let $f,g:(X,0)\rightarrow(\mathbb{C},0)$ be complex analytic function-germs such that $f$ has isolated singularity at the origin. Let $\mathcal{W}$ be the Whitney stratification of $X$ and $\mathcal{V}$ be the good stratification of $X$ induced by $f.$ Suppose that $\Sigma_{\mathcal{W}}g$ is one-dimensional and that $\Sigma_{\mathcal{W}}g\cap\{f=0\}=\{0\}.$

Suppose that $g$ is tractable at the origin with respect to $\mathcal{V}.$ By \cite{Santana},
\begin{eqnarray}
\mathcal{V}^{\prime}=\Big\{ V_{i}\setminus \Sigma_{\mathcal{W}} g, V_{i}\cap \Sigma_{\mathcal{W}} g, V_i\in\mathcal{V}\Big\} \cup\mathcal{V}^f\label{good stratification}
\end{eqnarray} 
\noindent is a good stratification of $X$ relative to $f,$ where  $\mathcal{V}^f$ denotes the collection of strata of $\mathcal{V}$ contained in $\{f=0\}$ and
\begin{eqnarray*}
\mathcal{V}^{\prime\prime}=\{V_i\setminus\{g=0\}, V_i\cap\{g=0\}\setminus\Sigma_{\mathcal{W}}g, V_i\cap\Sigma_{\mathcal{W}}g, V_i\in\mathcal{V}\}\cup\{0\},\label{ second good stratification}
\end{eqnarray*} 
\noindent is a good stratification of $X$ relative to $g.$ Let us denote by $\tilde{\mathcal{V}}$ the good stratification of $X$ induced by $\tilde{g}=g+f^N, N\gg1$.

Let $\alpha$ be a regular value of $g,$  $\alpha^{\prime}$ a regular value of $\tilde{g},$ $ 0<|\alpha|,|\alpha^{\prime}|\ll\epsilon\ll1,$ $n$ be the number of stratified Morse points of a Morsefication of $f|_{X\cap g^{-1}(\alpha)\cap B_{\epsilon}}$ appearing on \linebreak$X_{reg}\cap g^{-1}(\alpha)\cap\{f\neq0\}\cap B_{\epsilon}$ and $\tilde{n}$ be the number of stratified Morse points of a Morsefication of $f|_{X\cap \tilde{g}^{-1}(\alpha^{\prime})\cap B_{\epsilon}}$ appearing on $X_{reg}\cap \tilde{g}^{-1}(\alpha^{\prime})\cap\{f\neq0\}\cap B_{\epsilon}.$ 

\begin{proposition}\label{comparing Brasselet number, n and n tilde}
Suppose that $g$ is tractable at the origin with respect to $\mathcal{V}.$ Then, \begin{center}
$B_{g,X}(0)-B_{\tilde{g},X}(0)=(-1)^{d-1}(n-\tilde{n}).$
\end{center}
\end{proposition}

\noindent\textbf{Proof.} By \cite{Santana}, \begin{center}
$B_{g,X}(0)-B_{f,X}(0)=(-1)^{d-1}(n-m)-\sum_{j=1}^{r}m_{f,b_j}(Eu_{X}(b_j)-B_{g,X\cap\{f=\delta\}}(b_j)),$
\end{center}
\noindent where $m$ is the number of stratified Morse points of a Morsification of $g|_{X\cap f^{-1}(\delta)\cap B_{\epsilon}}$ appearing on $X_{reg}\cap f^{-1}(\delta)\cap\{g\neq0\}\cap B_{\epsilon},$ for $0<|\delta|\ll\epsilon\ll 1.$

By Lemma \ref{g tilde is prepolar}, $\tilde{g}$ is prepolar at the origin with respect to $\mathcal{V}.$ So, by Corollary 6.5 of \cite{DG}, \begin{center}
    $B_{\tilde{g},X}(0)-B_{f,X}(0)=(-1)^{d-1}(\tilde{n}-\tilde{m}),$
\end{center}
\noindent where $\tilde{m}$ is the number of stratified Morse points of a Morsification of $\tilde{g}|_{X\cap f^{-1}(\delta)\cap B_{\epsilon}}$ appearing on $X_{reg}\cap f^{-1}(\delta)\cap\{\tilde{g}\neq0\}\cap B_{\epsilon}.$

Using Corollary \ref{Comparing m and m tilde} and Theorem \ref{Euler obstruction of a function formula}, we have the formula.\fim

\begin{lemma}\label{intersection multiplicity 2}
Suppose that $g$ is tractable at the origin with respect to $\mathcal{V}$ relative to $f.$ If $N\gg1$ is bigger than the maximum gap ratio of all components of the symmetric relative polar curve $\tilde{\Gamma}_{f,g}(\mathcal{V})$ and such that Proposition \ref{g tilde has isolated singularity} is satisfied, then \begin{eqnarray*}
    \left([\tilde{\Gamma}_{f,g}(\mathcal{V})].[V(g)]\right)_0=\left([\tilde{\Gamma}_{f,g}(\mathcal{V})].[V(\tilde{g})]\right)_0.
\end{eqnarray*}
\end{lemma}
\noindent\textbf{Proof.} Since $g$ is tractable at the origin with respect to $\mathcal{V}, \tilde{\Gamma}_{f,g}(\mathcal{V})$ is a curve. Let us write $[\tilde{\Gamma}_{f,g}(\mathcal{V})]=\sum_{v}m_v[v],$ where each component $v$ of $\tilde{\Gamma}_{f,g}(\mathcal{V})$ is a reduced irreducible curve at the origin. Let $\alpha_v(t)$ be a parametrization  of $v$ such that $\alpha_v(0)=0.$ By page 974 \cite{Ms1}, each component $v$ intersects $V(g-g(p))$ at a point $p\in v, p\neq 0$, sufficiently close to the origin and such that $g(p)\neq 0.$ So, \begin{equation*}
    codim_X\{0\}=codim_X V(g)+codim_X v .
\end{equation*}

Also, each component (reduced irreducible curve at the origin) $\tilde{v}$ of $\tilde{\Gamma}_{f,\tilde{g}}(\mathcal{V})$ intersects $V(\tilde{g}-\tilde{g}(p))$ at such point $p\in \tilde{v}, p\neq 0$ and $\tilde{g}(p)\neq0.$ But since $\tilde{\Gamma}_{f,\tilde{g}}(\mathcal{V})=\tilde{\Gamma}_{f,g}(\mathcal{V})\cup\Sigma_{\mathcal{W}} g,$ we also have that $v$ intersects $V(\tilde{g}-\tilde{g}(p))$ at the point $p,$  so \begin{equation*}
    codim_X\{0\}=codim_X V(\tilde{g})+codim_X v.
\end{equation*}

Therefore, by A.9 of \cite{massey2003numerical} , \begin{eqnarray}
\left([v].[V(g)]\right)_0&=&mult_t g(\alpha_v(t))\nonumber\\
\left([v].[V(\tilde{g})]\right)_0&=&min\{mult_t g(\alpha_v(t)),mult_t f^N(\alpha_v(t))\}\nonumber
\end{eqnarray}
Now, \begin{center} $mult_t f^N(\alpha_v(t))=N \left([v].[V(f)]\right)_0$ and $mult_t g(\alpha_v(t))= \left([v].[V(g)]\right)_0.$
\end{center}
The gap ratio of $v$ at the origin for $g$ with respect to $f$ is the ratio of intersection numbers $\frac{\left([v].[V(g)]\right)_0}{\left([v].[V(f)]\right)_0}.$ So, if $N>\frac{\left([v].[V(g)]\right)_0}{\left([v].[V(f)]\right)_0},$ then $mult_t f^N(\alpha_v(t))>mult_t g(\alpha_v(t)).$

Making the same procedure over each component $v$ of $\tilde{\Gamma}_{f,g}(\mathcal{V})$ and using that $N$ is bigger then the maximum gap ratio of all components $v$ of $\tilde{\Gamma}_{f,g}(\mathcal{V})$ and such that Proposition \ref{g tilde has isolated singularity} is satisfied, we conclude that \begin{eqnarray*}
    \left([\tilde{\Gamma}_{f,g}(\mathcal{V})].[V(g)]\right)_0=\left([\tilde{\Gamma}_{f,g}(\mathcal{V})].[V(\tilde{g})]\right)_0.
\end{eqnarray*}
\fim

Our next goal is give another proof for the Lê-Iomdin formula for the Brasselet number. For that we need to compare $n$ and $\tilde{n}.$ 

\begin{lemma}\label{Comparing n and n tilde}
If $N$ is bigger than the maximum gap ratio of all components of the symmetric relative polar curve $\tilde{\Gamma}_{f,g}(\mathcal{V})$ and such that Proposition \ref{g tilde has isolated singularity} is satisfied, if  $0<|\alpha|,|\alpha^{\prime}|\ll\epsilon\ll1,$ then
\begin{eqnarray*}
    \tilde{n}=n+ (-1)^{d-1}N\sum_{j=1}^{r}m_{f,b_j} Eu_{f,X\cap \tilde{g}^{-1}(\alpha^{\prime})}(b_j).
\end{eqnarray*}
\end{lemma}
\noindent\textbf{Proof.} We start describing the critical points of $f|_{g^{-1}(\alpha)\cap B_{\epsilon}}.$ We have
\vspace{-0,1cm}
\begin{eqnarray}
x\in\Sigma f|_{g^{-1}(\alpha)\cap B_{\epsilon}} &\Leftrightarrow& x\in g^{-1}(\alpha)\cap B_{\epsilon} \ and \ rk(d_xg,d_xf)\leq1\nonumber\\
 &\Leftrightarrow & x\in g^{-1}(\alpha)\cap B_{\epsilon} \ and \ \big(d_xg=0\big) \ or \ \big(d_xf=0\big) \ or \  \big(d_xg=\lambda d_xf, \lambda\neq0\big).\nonumber
\end{eqnarray}
Since $f$ has isolated singularity at the origin and, by Proposition 1.3 of \cite{Ms1}, $\Sigma_{\mathcal{W}} g\subset\{g=0\},$ we have that $\Sigma f|_{g^{-1}(\alpha)\cap B_{\epsilon}}=g^{-1}(\alpha)\cap B_{\epsilon}\cap\tilde{\Gamma}_{f,g}(\mathcal{V}).$ Therefore, $n$ counts the number of Morse points of a Morsification of $f|_{g^{-1}(\alpha)\cap B_{\epsilon}}$ coming from $g^{-1}(\alpha)\cap B_{\epsilon}\cap\tilde{\Gamma}_{f,g}(\mathcal{V}).$ 

Now, let us describe $\Sigma f|_{\tilde{g}^{-1}(\alpha^{\prime})\cap B_{\epsilon}}.$
\begin{eqnarray}
x\in\Sigma f|_{\tilde{g}^{-1}(\alpha^{\prime})\cap B_{\epsilon}} &\Leftrightarrow& x\in \tilde{g}^{-1}(\alpha^{\prime})\cap B_{\epsilon} \ and \ rk(d_x\tilde{g},d_xf)\leq1\nonumber\\
&\Leftrightarrow & x\in \tilde{g}^{-1}(\alpha^{\prime})\cap B_{\epsilon} \ and \ \big(d_x\tilde{g}=0\big) \ or \ \big(d_xf=0\big) \ or \  \big(d_x\tilde{g}=\lambda^{\prime} d_xf,\lambda^{\prime}\neq 0\big).\nonumber
\end{eqnarray}

Since $f$ and $\tilde{g}$ has isolated singularity at the origin, we have that \begin{equation*}
    \Sigma f|_{\tilde{g}^{-1}(\alpha^{\prime})\cap B_{\epsilon}}=\tilde{g}^{-1}(\alpha^{\prime})\cap B_{\epsilon}\cap\tilde{\Gamma}_{f,\tilde{g}}(\mathcal{V}).
\end{equation*}

Since $\tilde{\Gamma}_{f,\tilde{g}}(\mathcal{V})=\tilde{\Gamma}_{f,g}(\mathcal{V})\cup\Sigma_{\mathcal{W}} g,$  \begin{equation*}
    \Sigma f|_{\tilde{g}^{-1}(\alpha^{\prime})\cap B_{\epsilon}}=\big(\Sigma_{\mathcal{W}} g\cap\tilde{g}^{-1}(\alpha^{\prime})\cap B_{\epsilon}\big)\cup\big(\tilde{\Gamma}_{f,g}(\mathcal{V})\cap\tilde{g}^{-1}(\alpha^{\prime})\cap B_{\epsilon}\big).
\end{equation*}

Notice that, since $\Sigma_{\mathcal{W}} g\cap\{f=0\}=\{0\}, \Sigma_{\mathcal{W}} g\cap\tilde{g}^{-1}(\alpha^{\prime})\cap B_{\epsilon}\subset\{f\neq0\}.$ Also, by definition, $\tilde{\Gamma}_{f,g}(\mathcal{V})\setminus\{0\}\subset\{f\neq0\}$ Therefore, $\tilde{n}$ counts the number of Morse points of a Morsification of $f|_{\tilde{g}^{-1}(\alpha^{\prime})\cap B_{\epsilon}}$ coming from $\tilde{g}^{-1}(\alpha^{\prime})\cap B_{\epsilon}\cap\Sigma_{\mathcal{W}} g\cap \{f\neq0\}\cap \{g=0\}$ and from \linebreak$\tilde{g}^{-1}(\alpha^{\prime})\cap B_{\epsilon}\cap\tilde{\Gamma}_{f,g}(\mathcal{V})\cap \{f\neq0\}\cap \{g\neq0\}.$ 

By Lemma \ref{intersection multiplicity 2}, the number of Morse points of a Morsification of $f|_{\tilde{g}^{-1}(\alpha^{\prime})\cap B_{\epsilon}}$ appearing on $\tilde{g}^{-1}(\alpha^{\prime})\cap B_{\epsilon}\cap\tilde{\Gamma}_{f,g}(\mathcal{V})\cap \{f\neq0\}\cap \{g\neq0\}$ is precisely $n.$ Let us describe the number of Morse points of a Morsification of $f|_{\tilde{g}^{-1}(\alpha^{\prime})\cap B_{\epsilon}}$ appearing on \linebreak$\tilde{g}^{-1}(\alpha^{\prime})\cap B_{\epsilon}\cap\Sigma_{\mathcal{W}} g\cap \{f\neq0\}\cap \{g=0\}.$
Using that $\Sigma_{\mathcal{W}} g\subset\{g=0\},$
\begin{eqnarray}
x\in\tilde{g}^{-1}(\alpha^{\prime})\cap B_{\epsilon}\cap\Sigma_{\mathcal{W}} g &\Leftrightarrow& \tilde{g}(x)=\alpha^{\prime} \ and \ d_xg=0\nonumber\\
&\Leftrightarrow& g(x)+f(x)^N=\alpha^{\prime} \ and \ d_xg=0\nonumber\\
&\Leftrightarrow& f(x)^N=\alpha^{\prime} \ and \ d_xg=0\nonumber\\
&\Leftrightarrow& f(x)\in\{\alpha_0,\ldots,\alpha_{N-1}\} \ and \ d_xg=0,\nonumber
\end{eqnarray}
\noindent where $\{\alpha_0,\ldots,\alpha_{N-1}\}$ are the $N$-th roots of $\alpha^{\prime}.$ 
Therefore, \begin{equation*}
    \tilde{g}^{-1}(\alpha^{\prime})\cap B_{\epsilon}\cap\Sigma_{\mathcal{W}} g=\bigcup_{i=0}^{N-1}f^{-1}(\alpha_i)\cap B_{\epsilon}\cap\Sigma_{\mathcal{W}} g.
\end{equation*}

Since $\Sigma_{\mathcal{W}}g$ is one-dimensional, $f^{-1}(\alpha_i)\cap\Sigma_{\mathcal{W}}g$ is a finite set of critical points of $f|_{\tilde{g}^{-1}(\alpha^{\prime})\cap B_{\epsilon}}.$
Since $\tilde{\Gamma}_{f,\tilde{g}}(\mathcal{V})=\Sigma_{\mathcal{W}} g\cup\tilde{\Gamma}_{f,g}(\mathcal{V}),$ each branch $b_j$ of $\Sigma_{\mathcal{W}} g$ is a component of $\tilde{\Gamma}_{f,\tilde{g}}(\mathcal{V}).$ If $V_{i(j)}$ is the stratum of $\mathcal{V}^{\prime\prime}$ containing $b_j,$ then $f|_{V_{i(j)}\cap \tilde{g}^{-1}(\alpha^{\prime})}$ has isolated singularity at each point $x_{l}\in b_j\cap f^{-1}(\alpha_{i})\cap \tilde{g}^{-1}(\alpha^{\prime}), j\in\{1,\ldots,r\}$ and $l\in\{i_1,\ldots,i_{k(j)}\}$ (page 974, \cite{Ms1}). Using Proposition \ref{Eu_f and Morse points}, we can count the number $n_{l}$ of Morse points of a Morsification of $f|_{\tilde{g}^{-1}(\alpha^{\prime})\cap B_{\epsilon}}$ in a neighborhood of each $x_{l},$ \begin{equation*}
    Eu_{f,X\cap \tilde{g}^{-1}(\alpha^{\prime})}(x_{l})=(-1)^{d-1}n_{l}.
\end{equation*}

Since the Euler obstruction of a function is constant on each branch $b_j$, so is the Euler obstruction of a function and we can denote $Eu_{f,X\cap \tilde{g}^{-1}(\alpha^{\prime})}(x_{l})$ by $Eu_{f,X\cap \tilde{g}^{-1}(\alpha^{\prime})}(b_j),$ for all $x_{l}\in b_j\cap f^{-1}(\alpha_{i})\cap \tilde{g}^{-1}(\alpha^{\prime}).$
Therefore, if $b_j\cap f^{-1}(\alpha_{i})\cap \tilde{g}^{-1}(\alpha^{\prime})=\{x_{j_1},\ldots, x_{j_{m_{f,b_j}}}\},$ the number of Morse points of a Morsification of $f|_{\tilde{g}^{-1}(\alpha^{\prime})\cap B_{\epsilon}}$ appearing on $\big(X_{reg}\setminus\{\tilde{g}=0\}\big)\cap b_j\cap\{\tilde{g}=\alpha^{\prime}\}\cap B_{\epsilon}\cap \{f=\alpha_i\}$ is \begin{eqnarray}
    n_{j_1}+\cdots+n_{j_{m_{f,b_j}}}&=&(-1)^{d-1}{m_{f,b_j}} Eu_{f,X\cap \tilde{g}^{-1}(\alpha^{\prime})}(x_{l}).\nonumber
\end{eqnarray}
\noindent
Making the same analysis over each $\alpha_i\in\sqrt[N]{\alpha^{\prime}},$ the number of Morse points of a Morsification of $f|_{\tilde{g}^{-1}(\alpha^{\prime})\cap B_{\epsilon}}$ appearing in $X_{reg}\setminus\{\tilde{g}=0\}\cap\{g=0\}\cap\{\tilde{g}=\alpha^{\prime}\}\cap B_{\epsilon}$ is \begin{equation*}
    (-1)^{d-1}N\sum_{j=1}^r m_{f,b_j}  Eu_{f,X\cap \tilde{g}^{-1}(\alpha^{\prime})}(b_j).
\end{equation*}
\fim

\begin{theorem}\label{Le-Iomdin formula for the Brasselet number}
Suppose that $g$ is tractable at the origin with respect to $\mathcal{V}.$ If $0<|\alpha|,|\alpha^{\prime}|\ll\epsilon$ and $N$ is bigger than the maximum gap ratio of each component of the symmetric relative polar curve $\tilde{\Gamma}_{f,g}(\mathcal{V})$ and such that Proposition \ref{g tilde has isolated singularity} is satisfied, then \begin{eqnarray*}
    B_{\tilde{g},X}(0)=B_{g,X}(0)+ N\sum_{j=1}^{r}m_{f,b_j} Eu_{f,X\cap \tilde{g}^{-1}(\alpha^{\prime})}(b_j).
\end{eqnarray*}
\end{theorem}
\noindent\textbf{Proof.} It follows by Proposition \ref{comparing Brasselet number, n and n tilde} and Lemma \ref{Comparing n and n tilde}.\fim

This formula gives a way to compare the numerical data associated to the generalized Milnor fibre of a function $g$ with a one-dimensional singular locus and and to the generalized Milnor fibre of the deformation $\tilde{g}=g+f^N,$ for $N\gg1$ sufficiently large. This is what Lê \cite{le1980ensembles} and Iomdin \cite{iomdin1974complex} have done in the case where $g$ is defined over a complete intersection in $\mathbb{C}^n,$ $g$ has a one-dimensional critical locus and $f$ is a generic linear form over $\mathbb{C}^n.$ Therefore, Theorem \ref{Le-Iomdin formula for the Brasselet number} generalizes this Lê-Iomdin formula.

For $X=\mathbb{C}^n,$ let us consider $\mathcal{W}=\{\mathbb{C}^n\setminus\{0\},\{0\}\}$ the Whitney stratification of $\mathbb{C}^n.$ If $f$ has isolated singularity at the origin, the good stratification $\mathcal{V}$ of $\mathbb{C}^n$ induced by $f$ is given by $\mathcal{V}=\{\mathbb{C}^n\setminus\{f=0\}, \{f=0\}\setminus\{0\},\{0\}\}.$ 
\begin{corollary}\label{Le-Iomdin formula as a consquence}
Suppose that $g$ is tractable at the origin with respect to $\mathcal{V}$ relative to $f.$ If $\alpha$ and $\alpha^{\prime}$ are regular values of $g$ and $\tilde{g},$ respectively, with $0<|\alpha|,|\alpha^{\prime}|\ll\epsilon,$ then \begin{eqnarray*}
   \chi(\tilde{g}^{-1}(\alpha^{\prime})\cap B_{\epsilon})=\chi(g^{-1}(\alpha)\cap B_{\epsilon})+(-1)^{n-1}N\sum_{j=1}^{r}m_{f,b_j} \mu(g|_{f^{-1}(\delta_{j_i})},b_j),
\end{eqnarray*}
\noindent where $\mu(g|_{f^{-1}(\delta_{j_i})},b_j)$ denotes the Milnor number of $g|_{X\cap f^{-1}(\delta_{j_i})\cap B_{\epsilon}}$ at a point $x_{j_i}$ of the branch $b_j,$ with $f(x_{j_i})=\delta_{j_i}.$ 
\end{corollary}

\noindent\textbf{Proof.}  By \cite{Santana}, $\mathcal{V}^{\prime}=\{\mathbb{C}^n\setminus\{f=0\}\cup\Sigma_{\mathcal{W}}g,\{f=0\}\setminus\{0\},\Sigma_{\mathcal{W}} g,\{0\}\}$ is a good stratification of $\mathbb{C}^n$ relative to $f.$ Also by \cite{Santana}, $\mathcal{V}^{\prime\prime},$ given by \begin{equation*}
    \{\mathbb{C}^n\setminus\{f=0\}\cup\{g=0\},\{f=0\}\setminus\{g=0\},\{g=0\}\setminus\{f=0\}\cup\Sigma_{\mathcal{W}}g, \{f=0\}\cap\{g=0\}\setminus\Sigma_{\mathcal{W}}g, \Sigma_{\mathcal{W}}g,\{0\}\},
\end{equation*}
\noindent is a good stratification of $\mathbb{C}^n$ relative to $g.$

By definition of the Brasselet number, if $0<|\alpha|\ll\epsilon\ll1,$ \begin{eqnarray}
    B_{g,X}(0)&=&\sum_{V_i\in\mathcal{V}^{\prime\prime}}\chi(V_i\cap g^{-1}(\alpha)\cap B_{\epsilon})Eu_{\mathbb{C}^n}(V_i)\nonumber\\
    &=&\chi\big((\mathbb{C}^n\setminus\{f=0\}\cup\{g=0\})\cap g^{-1}(\alpha)\cap B_{\epsilon}\big)Eu_{\mathbb{C}^n}(\mathbb{C}^n\setminus\{f=0\}\cup\{g=0\})\nonumber\\
    &+&\chi\big((\{f=0\}\setminus\{g=0\})\cap g^{-1}(\alpha)\cap B_{\epsilon}\big)Eu_{\mathbb{C}^n}(\{f=0\}\setminus\{g=0\})\nonumber\\
    &=&\chi\big((\mathbb{C}^n\setminus\{g=0\})\cap g^{-1}(\alpha)\cap B_{\epsilon}\big)\nonumber\\
    &=&\chi(g^{-1}(\alpha)\cap B_{\epsilon}).\nonumber
\end{eqnarray}

The good stratification of $\mathbb{C}^n$ induced by $\tilde{g}$ is $\tilde{\mathcal{V}}=\{\{\tilde{g}=0\},\mathbb{C}^n\setminus\{\tilde{g}=0\},\{0\}\}$ and  then, if $0<|\alpha^{\prime}|\ll\epsilon\ll1,$ \begin{equation*}
    B_{\tilde{g},X}(0)=\chi(\mathbb{C}^n\setminus\{\tilde{g}=0\}\cap g^{-1}(\alpha)\cap B_{\epsilon})Eu_{\mathbb{C}^n}(\mathbb{C}^n\setminus\{0\})=\chi(\tilde{g}^{-1}(\alpha^{\prime})\cap B_{\epsilon}).
\end{equation*} 

Since $f|_{\tilde{g}^{-1}(\alpha^{\prime})\cap B_{\epsilon}}$ is defined over $\mathbb{C}^n$ and has isolated singularity at each $x_{j_i}\in b_j$, \linebreak considering a small ball $B_{\epsilon}(x_{j_i})$ with radius $\epsilon$ and center at $x_{j_i},$ for $0<|\delta|\ll\epsilon\ll1,$ \begin{eqnarray}
    Eu_{f,\tilde{g}^{-1}(\alpha^{\prime})}(x_{j_i})&=&(-1)^{n-1}\mu\big(f|_{\tilde{g}^{-1}(\alpha^{\prime})},x_{j_i}\big)\nonumber\\
    &=&(-1)^{n-1}(-1)^{n-1}\big[\chi\big((f|_{\tilde{g}^{-1}(\alpha^{\prime})})^{-1}(\delta)\cap B_{\epsilon}(x_{j_i})\big)-1\big]\nonumber\\
&=&\chi\big(f^{-1}(\delta_{j_i}-\delta)\cap \tilde{g}^{-1}(\alpha^{\prime})\cap B_{\epsilon}(x_{j_i})\big)-1, \ f(x_{j_i})=\delta_{j_i}\nonumber\\
&\stackrel{*}{=}&\chi\big(f^{-1}(\delta_{j_i})\cap \tilde{g}^{-1}(\alpha^{\prime}-\delta)\cap B_{\epsilon}(x_{j_i})\big)-1 \nonumber\\
&=&\chi\big(f^{-1}(\delta_{j_i})\cap g^{-1}(\alpha^{\prime}-\delta_{j_i}^N-\delta)\cap B_{\epsilon}(x_{j_i})\big)-1, g(x_{j_i})=\alpha^{\prime}-\delta_{j_i}^N\nonumber\\
&=&\chi\big((g|_{f^{-1}(\delta_{j_i})})^{-1}(\delta)\cap B_{\epsilon}(x_{j_i})\big)-1\nonumber\\
&=&(-1)^{n-1}\mu\big(g|_{f^{-1}(\delta_{j_i})},x_{j_i}\big),\nonumber
\end{eqnarray}
\noindent where the equality $(*)$ is justified by Proposition 6.2 of \cite{DG}. 
Therefore, applying Theorem \ref{Le-Iomdin formula for the Brasselet number}, we obtain \begin{eqnarray*}
   \chi(\tilde{g}^{-1}(\alpha^{\prime})\cap B_{\epsilon})=\chi(g^{-1}(\alpha)\cap B_{\epsilon})+(-1)^{n-1}N\sum_{j=1}^{r}m_{f,b_j} \mu(g|_{f^{-1}(\delta_{j_i})},b_j).
\end{eqnarray*} \fim

Another consequence of Theorem \ref{Le-Iomdin formula for the Brasselet number} is a different proof for the Lê-Iomdin formula proved by Massey in \cite{massey2003numerical} in the case of a function with a one-dimensional singular locus. For that we will need the definition of the Lê-numbers.
 We present here the case for functions defined over a nonsingular subspace of $\mathbb{C}^n$, and we recommend Part I of \cite{massey2003numerical} for the general case. Let $h:(U,0)\subseteq(\mathbb{C}^n,0)\rightarrow(\mathbb{C},0)$ be an analytic function such that its critical locus $\Sigma h$ is a $s$-dimensional set. For $0\leq k\leq n,$ the \textbf{$k$-th relative polar variety} $\Gamma^k_{h,z}$ of $h$ with respect to $z$ is the scheme $V\left(\frac{\partial h}{\partial z_k},\ldots,\frac{\partial h}{\partial z_n}\right)/\Sigma h,$ where $z=(z_1,\ldots,z_n)$ are fixed local coordinates and the \textbf{$k$-th polar cycle} of $h$ with respect to $z$ is the analytic cycle $[\Gamma^{k}_{h,z}].$ The \textbf{$k$-th Lê cycle} $[\Lambda^k_{h,z}]$ of $h$ with respect to $z$ is the difference of cycles $\left[\Gamma^{k+1}_{h,z}\cap V(\frac{\partial h}{\partial z_k})\right]-[\Gamma^{k}_{h,z}].$

\begin{definition}\label{le number definition}
The \textbf{$k$-th Lê number} of $h$ in $p$ with respect to $z, \ \lambda^k_{h,z},$ is the intersection number\begin{equation*}
    (\Lambda^k_{h,z}.V(z_0-p_0,\ldots,z_{k-1}-p_{k-1}))_p,
\end{equation*} provided this intersection is purely zero-dimensional at $p.$

If this intersection is not purely zero-dimensional, the $k$-th Lê number of $h$ at $p$ with respect to $z$ is said to be undefined. 
\end{definition}
\begin{corollary}\label{Massey Le-Iomdin formula as a consequence}
 Let $\mathcal{V}$ be the good stratification of an open set $(U,0)\subseteq(\mathbb{C}^{n+1},0)$ induced by a generic linear form $l$ defined over $\mathbb{C}^{n+1}.$ Let $N\geq2, \mathbf{z}=(z_0\ldots,z_n)$ be a linear choice of coordinates such that $\lambda^i_{g,\mathbf{z}}(0)$ is defined for $i=0,1,$ and $\mathbf{\tilde{z}}=(z_1\ldots,z_n,z_0)$ be the coordinates for $\tilde{g}=g+l^N$ such that $\lambda^0_{\tilde{g},\mathbf{\tilde{z}}}$ is defined. If $N$ is greater then the maximum gap ratio of each component of the symmetric relative polar curve $\tilde{\Gamma}_{f,g}(\mathcal{V})$ and such that Proposition \ref{g tilde has isolated singularity} is satisfied, then \begin{eqnarray*}
    \lambda^0_{\tilde{g},\mathbf{\tilde{z}}}(0)=\lambda^0_{g,\mathbf{z}}(0)+ (N-1)\lambda^1_{g,\textbf{z}}(0).
\end{eqnarray*}
\end{corollary}

\noindent\textbf{Proof.} By \cite{Santana}, $g$ is tractable at te origin with respect to the good stratification $\mathcal{V}$ induced by $l.$ Without loss of generality, we can suppose that $l=z_0.$ Let $F_{g,0}$ be the Milnor fibre of $g$ at the origin and $F_{\tilde{g},0}$ be the Milnor fibre of $\tilde{g}$ at the origin. Since $g$ has a one-dimensional critical set, the possibly nonzero Lê numbers are $\lambda^0_{g,\mathbf{z}}(0)$ and $\lambda^1_{g,\mathbf{z}}(0)$ and, since $\tilde{g}$ has isolated singularity at the origin, the only possibly nonzero Lê number is $\lambda^0_{\tilde{g},\mathbf{\tilde{z}}}(0).$ By Theorem 4.3 of \cite{massey1988le}, \begin{equation*}
    \chi(F_{g,0})=1+(-1)^n\lambda^0_{g,\mathbf{z}}(0)+(-1)^{n-1}\lambda^1_{g,\mathbf{z}}(0)
\end{equation*}
and 
\begin{equation*}
    \chi(F_{\tilde{g},0})=1+(-1)^n\lambda^0_{\tilde{g},\mathbf{\tilde{z}}}(0).
\end{equation*}

In \cite{massey2003numerical}, on page 49, Massey remarked that for $0<|\delta|\ll\epsilon\ll1,$
\begin{equation*}
    \lambda^1_{g,\mathbf{z}}(0)=\sum_{j=1}^{r}m_{b_j} \mu(g|_{l^{-1}(\delta)},b_j).
\end{equation*}

Therefore, the formula holds by Corollary \ref{Le-Iomdin formula as a consquence}.\fim

\section{Applications to generic linear forms}

\hspace{0,5cm} Let $g:(X,0)\rightarrow(\mathbb{C},0)$ be a complex analytic function-germ and $l$ be a generic linear form in $\mathbb{C}^n.$ Let $\mathcal{W}=\{\{0\},W_1,\ldots, W_q\}$ be a Whitney stratification of $X$ and $\mathcal{V}$ be the good stratification of $X$ induced by $l.$ Suppose that $\Sigma_{\mathcal{W}}g$ is one-dimensional.

Let $\mathcal{V}^{\prime}$ be the good stratification of $X$ relative to $l,$ $\mathcal{V}^{\prime\prime}$ be the good stratification of $X$ relative to $g$ and $\tilde{\mathcal{V}}$ be the good stratification of $X$ induced by $\tilde{g}=g+l^N, N\gg1$, taken as in Section 4. 

Let $\alpha$ be a regular value of $g,$  $\alpha^{\prime}$ a regular value of $\tilde{g},$ $ 0<|\alpha|,|\alpha^{\prime}|\ll\epsilon\ll1,$ $n$ be the number of stratified Morse points of a Morsification of $l|_{X\cap g^{-1}(\alpha)\cap B_{\epsilon}}$ appearing on \linebreak$X_{reg}\cap g^{-1}(\alpha)\cap\{l\neq0\}\cap B_{\epsilon},$ $n_i$ be the number of stratified Morse points of a Morsification of $l|_{W_i\setminus(\{g=0\}\cup\{l=0\})\cap g^{-1}(\alpha)\cap B_{\epsilon}}$ appearing on $W_i\cap g^{-1}(\alpha)\cap\{l\neq0\}\cap B_{\epsilon},$ $\tilde{n}$ be the number of stratified Morse points of a Morsification of $l|_{X\cap \tilde{g}^{-1}(\alpha^{\prime})\cap B_{\epsilon}}$ appearing on $X_{reg}\cap \tilde{g}^{-1}(\alpha^{\prime})\cap\{l\neq0\}\cap B_{\epsilon}$ and $\tilde{n}_i$ be the number of stratified Morse points of a Morsification of $l|_{W_i\setminus\{\tilde{g}=0\}\cap \tilde{g}^{-1}(\alpha^{\prime})\cap B_{\epsilon}}$ appearing on $W_i\cap \tilde{g}^{-1}(\alpha^{\prime})\cap\{l\neq0\}\cap B_{\epsilon},$  for each $W_i\in\mathcal{W}.$

As before, we write $\Sigma_{\mathcal{W}} g$ as a union of branches $b_1\cup\ldots\cup b_r$ and we suppose that $\{l=\delta\}\cap b_j=\{x_{i_1},\ldots,x_{i_{k(j)}}\}.$ For each $t\in\{i_1,\ldots,i_{k(j)}\},$ let $D_{x_t}$ be the closed ball with center at $x_{t}$ and radius $r_{t}, 0<|\alpha|,|\alpha^{\prime}|\ll|\delta|\ll r_{t}\ll\epsilon\ll1,$ sufficiently small for the balls $D_{x_t}$ be pairwise disjoint and the union of balls $D_j=D_{x_{i_1}}\cup\ldots\cup D_{x_{i_{k(j)}}}$ be contained in $B_{\epsilon}$ and $\epsilon$ is sufficiently small such that the local Euler obstruction of $X$ at a point of $b_j\cap B_{\epsilon}$ is constant.

In \cite{tibuar1998embedding}, Tib\u{a}r gave a bouquet decomposition to for the Milnor fibre of $\tilde{g}$ in terms of the Milnor fibre of $g.$ Let us denote by $F_g$ the local Milnor fibre of $g$ at the origin, $F_{\tilde{g}}$ the local Milnor fibre of $\tilde{g}$ at the origin and $F_j$ is the local Milnor fibre of $g|_{\{l=\delta\}}$ at a point of the branch $b_j.$ Then there is a homotopy equivalence \begin{center}
$F_{\tilde{g}}\stackrel{ht}{\simeq} (F_g\cup E)\bigvee_{j=1}^r\bigvee_{M_j}S(F_j),$
\end{center}
\noindent where $\bigvee$ denotes the wedge sum of topological spaces, $M_j=Nm_{b_j}-1, S(F_j)$ denotes the topological suspension over $F_j, E:=\cup_{j=1}^rCone(F_j)$ and $F_g\cup E$ is the attaching to $F_g$ of one cone over $F_j\subset F_g$ for each $j\in\{1,\ldots, r\}.$  As a consequence of this theorem, Tib\u{a}r proved a Lê-Iomdin formula for the Euler characteristic of these Milnor fibres. 





  

In the following, we present a new proof for this formula using our previous results. 
 
\begin{proposition}\label{proposition Nicolas-Tibar}
Suppose that $g$ is tractable at the origin with respect to $\mathcal{V}.$ If\linebreak $0<|\alpha|,|\alpha^{\prime}|\ll|\delta|\ll\epsilon\ll1,$ then
\begin{eqnarray*}
\chi(X\cap\tilde{g}^{-1}(\alpha^{\prime})\cap 
    B_{\epsilon})-\chi(X\cap g^{-1}(\alpha)\cap B_{\epsilon})
    &=&N\sum_{j=1}^{r}m_{b_j}\big(1-\chi (F_j)\big),
\end{eqnarray*}
\noindent where $F_j=X\cap g^{-1}(\alpha)\cap H_j\cap D_{x_{t}}$ is the local Milnor fibre of $g|_{\{l=\delta\}}$ at a point of the branch $b_j$ and $H_j$ denotes the generic hyperplane $l^{-1}(\delta)$ passing through $x_t\in b_j,$ for $t\in\{i_1,\ldots,i_{k(j)}\}.$
\end{proposition}
\noindent\textbf{Proof.} For a stratum $V_i=W_i\setminus(\{g=0\}\cup\{l=0\})$ in $\mathcal{V}^{\prime\prime},$  $W_i\in\mathcal{W}$, let $N_i$ be a normal slice to $V_i$ at $x_t\in b_j,$ for $t\in\{i_1,\ldots,i_{k(j)}\}$ and $D_{x_t}$ a closed ball of radius $r_t$ centered at $x_t$. Considering the constructible function $\textbf{1}_X,$ the normal Morse index along $V_i$ is given by  \begin{eqnarray*}
     \eta(V_i,\textbf{1}_X)&=&\chi(W_i\setminus(\{g=0\}\cup\{l=0\})\cap N_i\cap D_{x_t})\\&-&\chi(W_i\setminus(\{g=0\}\cup\{l=0\})\cap N_i\cap\{g=\alpha\}\cap D_{x_t})\\
    &=&\chi(W_i\cap N_i\cap D_{x_t})-\chi(W_i\cap N_i\cap\{g=\alpha\}\cap D_{x_t})\\
    &=&1-\chi(l_{W_i}).
\end{eqnarray*}

For a stratum $\tilde{V}_i=W_i\setminus(\{\tilde{g}=0\}\in\tilde{\mathcal{V}},$  $W_i\in\mathcal{W}$, let $\tilde{N}_i$ be a normal slice to $\tilde{V}_i$ at $x_t\in b_j,$ for $t\in\{i_1,\ldots,i_{k(j)}\}$. Considering the constructible function $\textbf{1}_X,$ the normal Morse index along $\tilde{V}_i$ is given by  \begin{eqnarray*}
     \eta(\tilde{V}_i,\textbf{1}_X)&=&\chi((W_i\setminus\{\tilde{g}=0\})\cap \tilde{N}_i\cap D_{x_t})-\chi((W_i\setminus\{\tilde{g}=0\})\cap \tilde{N}_i\cap\{\tilde{g}=\alpha^{\prime}\}\cap D_{x_t})\\
    &=&\chi(W_i\cap \tilde{N}_i\cap D_{x_t})-\chi(W_i\cap \tilde{N}_i\cap\{\tilde{g}=\alpha^{\prime}\}\cap D_{x_t})\\
    &=&1-\chi(l_{W_i}).
\end{eqnarray*}

Then applying Theorem 4.2 of \cite{DG} for $\textbf{1}_X$, we obtain that \begin{center}
    $\chi(X\cap\tilde{g}^{-1}(\alpha^{\prime})\cap 
    B_{\epsilon})-\chi(X\cap\tilde{g}^{-1}(\alpha^{\prime})\cap l^{-1}(0)\cap B_{\epsilon})=\sum_{i=1}^{q}(-1)^{d_i-1}\tilde{n}_i(1-\chi(l_{W_i}))$
\end{center}
 and that 
 \begin{center}
    
      $\chi(X\cap g^{-1}(\alpha)\cap B_{\epsilon})-\chi(X\cap g^{-1}(\alpha)\cap l^{-1}(0)\cap B_{\epsilon})=\sum_{i=1}^{q}(-1)^{d_i-1}n_i(1-\chi(l_{W_i})),$
 \end{center}
 \noindent where $d_i=\dim W_i.$
 
Therefore, since $\chi(X\cap\tilde{g}^{-1}(\alpha^{\prime})\cap l^{-1}(0)\cap B_{\epsilon})=\chi(X\cap g^{-1}(\alpha)\cap l^{-1}(0)\cap B_{\epsilon}),$

\begin{eqnarray*}
\chi(X\cap\tilde{g}^{-1}(\alpha^{\prime})\cap 
    B_{\epsilon})-\chi(X\cap g^{-1}(\alpha)\cap B_{\epsilon})=\sum_{i=1}^{q}(-1)^{d_i-1}(\tilde{n}_i-n_i)(1-\chi(l_{W_i})).
\end{eqnarray*}

Applying Lemma \ref{Comparing n and n tilde} and Corollary \ref{comparing euler obstruction of funcions}, we obtain, for each $i,$ \begin{eqnarray*}
    \tilde{n}_i&=&n_i+ (-1)^{d_i-1}N\sum_{j=1}^{r}m_{b_j} Eu_{l,\overline{W_i}\cap \tilde{g}^{-1}(\alpha^{\prime})}(b_j)\\
    &=& n_i+ (-1)^{d_i-1}N\sum_{j=1}^{r}m_{b_j} Eu_{g,\overline{W_i}\cap H_j}(b_j),
\end{eqnarray*} 
\noindent where $H_j$ denotes the generic hyperplane $l^{-1}(\delta)$ passing through $x_t\in b_j,$ for $t\in\{i_1,\ldots,i_{k(j)}\}.$

Hence \begin{eqnarray*}
\chi(X\cap\tilde{g}^{-1}(\alpha^{\prime})\cap 
    B_{\epsilon})-\chi(X\cap g^{-1}(\alpha)\cap B_{\epsilon})&=&N\sum_{i=1}^{q}\bigg(\sum_{j=1}^{r}m_{b_j} Eu_{g,\overline{W_i}\cap H_j}(b_j)\bigg)\big(1-\chi(l_{W_i})\big)\\
    &=&N\sum_{j=1}^{r}m_{b_j}\Big(1-\chi\big(X\cap g^{-1}(\alpha)\cap H_j\cap D_{x_{t}}\big)\Big)\\
    &=&N\sum_{j=1}^{r}m_{b_j}\big(1-\chi (F_j)\big),
\end{eqnarray*}
\noindent for $t\in\{i_1,\ldots,i_{k(j)}\}.$\fim





\vspace{2cm}
(Hellen Monção de Carvalho Santana) Universidade de São Paulo, Instituto de Ciências Matemáticas e de
Computação - USP, Avenida Trabalhador São-Carlense, 400 - Centro, São Carlos, Brazil.
\it{E-mail address}: hellenmcarvalho@hotmail.com

 \end{document}